\documentclass[12pt]{article}
\usepackage{hyperref}
\usepackage{amsmath}
\usepackage{amsfonts}
\usepackage{amssymb}
\begin{document}
%
%
\vskip 3mm

\noindent ON THE LIMIT BEHAVIOUR OF FINITE-SUPPORT BIVARIATE DISCRETE PRO\-BA\-BI\-LI\-TY DISTRI\-BU\-TIONS UNDER ITERATED PARTIAL SUMMATIONS 
\vskip 3mm

\vskip 5mm
\noindent L\'{i}via Le\v{s}\v{s}ov\'{a} and J\'{a}n Ma\v{c}utek

\noindent Department of Applied Mathematics and Statistics

\noindent Comenius University in Bratislava

\noindent Mlynsk\'{a} dolina, 842 48 Bratislava

\noindent livia.lessova@fmph.uniba.sk

\noindent jmacutek@yahoo.com

\vskip 3mm

\noindent ABSTRACT

Bivariate partial-sums discrete probability distributions are defined. The question of the existence of a limit distribution for iterated partial summations is solved for finite-support bivariate distributions which satisfy conditions under which the power method (known from matrix theory) can be used. An oscillating sequence of distributions, a phenomenon which has never been reported before, is presented.

\vskip 4mm

\noindent 1.   INTRODUCTION 

Let $\lbrace P_x^{(1)}\rbrace_{x=0}^{\infty}$ and $\lbrace P_x^*\rbrace_{x=0}^{\infty}$ be probability mass functions of two univariate discrete probability distributions defined on nonnegative integers. The distribution $\lbrace P_x^{(1)}\rbrace_{x=0}^{\infty}$ (the descendant distribution) is a partial-sums distribution created from $\lbrace P_x^*\rbrace_{x=0}^{\infty}$ (the parent distribution) if
\begin{equation}\label{eq:ps_univariate_definition}
P_x^{(1)}=c_1\sum_{j=x}^{\infty}g(j)P_j^*,\quad x=0,1,2,...,
\end{equation}
where $c_1$ is a normalization constant and $g(j)$ a real function. Several types of partial summations - for different choices of $g(j)$ - are mentioned in the comprehensive monograph by Johnson et al. (2005). An extensive survey of pairs of parents and descendants was provided by Wimmer and Altmann (2000). More detailed analyses (e.g. relations between probability generating functions of the parent and descendant distributions) can be found in Ma\v{c}utek (2003).

Partial summations from (\ref{eq:ps_univariate_definition}) can be applied iteratively. Take $\{P_{x}^{(1)}\}_{x=0}^{\infty}$, i.e. the descendant distribution from (\ref{eq:ps_univariate_definition}), as the parent, with function $g(j)$ remaining unaltered. We obtain the descendant of the second generation
\begin{equation*}
P_x^{(2)}=c_2\sum_{j=x}^{\infty}g(j)P^{(1)}_j,\quad x=0,1,2,...,
\end{equation*}
and, repeatedly applying the partial summation, for any $k\in \mathbb{N}$ the descendant of the $k$-th generation  
\begin{equation*}
P_x^{(k)}=c_k\sum_{j=x}^{\infty}g(j)P^{(k-1)}_j,\quad x=0,1,2,...,
\end{equation*}
$c_2$, $c_k$ being normalization constants.

The question whether the sequence of the descendant distributions has a limit was investigated by Ma\v{c}utek (2006) for a constant function $g(j)$. In this case, the answer is positive for a wide class of parent distributions, with the limit distribution being geometric. Ko\v{s}\v{c}ov\'a et al. (2018) presented a solution - albeit not a general one - of the problem if the parent distribution has a finite support.

In this paper we extend the result from Ko\v{s}\v{c}ov\'a et al. (2018) to bivariate discrete probability distributions.

\vskip 3mm

\noindent 2.  BIVARIATE PARTIAL-SUMS DISTRIBUTIONS

Research on partial-sums distributions is almost exclusively dedicated to univariate distributions (see Wimmer and Ma\v{c}utek (2012), and references therein). The only note on the bivariate (and $r$-variate) partial-sums distributions can be found in Kotz and Johnson (1991), who more or less restrict themselves to a suggestion to study multivariate cases.

Univariate partial-sums distributions from Section 1 can be naturally generalized to two dimensions as follows.

Let $\{P_{x,y}^*\}_{x,y=0}^{\infty}$ and $\{P_{x,y}^{(1)}\}_{x,y=0}^{\infty}$ be bivariate discrete distributions and let $g(x,y)$ be a real function. Then $\{P_{x,y}^{(1)}\}_{x,y=0}^{\infty}$ is the descendant of the parent $\{P_{x,y}^*\}_{x,y=0}^{\infty}$ if
\begin{equation}\label{eq:def_ps_bivariate}
P_{x,y}^{(1)}=c_1\sum_{i=x}^{\infty}\sum_{j=y}^{\infty}g(i,j)P_{i,j}^*.
\end{equation}
We will obtain the descendant of the $k$-th generation analogously to the univariate case, e.g. the $k$-th descendant is
\begin{equation*}
P_{x,y}^{(k)}=c_k\sum_{i=x}^{\infty}\sum_{j=y}^{\infty}g(i,j)P_{i,j}^{(k-1)}.
\end{equation*}   

We will show that if the parent distribution has a finite support of the size $m\times n$, the power method, which is a computational approach to finding matrix eigenvalues and eigenvectors, can in some cases be used to find the limit distribution.

\vskip 3mm

\noindent 3.  POWER METHOD AND ITS APPLICATION

The power method (see e.g. Golub and Van Loan (1996)) was suggested as a computational tool which enables, under certain conditions, to find an approximation of square matrix eigenvalues. The method can be applied to a diagonalizable matrix (i.e. a matrix which has linearly independent eigenvectors, or, equivalently, it is similar to a diagonal matrix) with a unique dominant eigenvalue (denote the eigenvalues $\lambda_1, \lambda_2, ..., \lambda_n$; there exists $k$ such that $\left|\lambda_k\right| > \left|\lambda_i\right|$, $i\neq k$). The eigenvector corresponding to the dominant eigenvalue is the dominant eigenvector.
	
If a matrix $A$ satisfies abovementioned conditions, then there exists a non-zero vector $x_0$ such that the sequence $\{A^kx_0\}_{k=1}^{\infty}$ converges to a multiple of the dominant eigenvector.

While the application of the power method is straightforward for univariate iterated partial sumations (see Ko\v{s}\v{c}ov\'a et al. (2018)), a bivariate distribution requires an additional step, namely, a vectorization of the probability matrix, which is, however, a standard operation in matrix theory (see e.g. Golub and Van Loan (1996)). Denote $\mathbb{P}^*$ the parent distribution, i.e.
\begin{equation*}
\mathbb{P}^*=
\left(\begin{matrix}.
P_{0,0}^* & P_{0,1}^* &\dots & P_{0,n-1}^* \\
P_{1,0}^* & P_{1,1}^* &\dots & P_{1,n-1}^* \\
 \vdots & \vdots & \ddots & \vdots \\
P_{m-1,0}^* & P_{m-1,1}^* &\dots & P_{m-1,n-1}^* 
\end{matrix}\right),
\end{equation*} 
the vectorization of $\mathbb{P}^*$ yields a vector of probabilities
\begin{equation*}
v(\mathbb{P}^*)=\left(P_{0,0}^*,\dots, P_{m-1,0}^*, P_{0,1}^*, \dots,  P_{m-1,1}^*,  P_{0,n-1}^*, \dots, P_{m-1,n-1}^* \right)^\mathbf{T}.
\end{equation*}
Now we will construct a matrix $\tilde{G}$ from the values of the function $g(i,j)$ from (\ref{eq:def_ps_bivariate}),
\begin{equation*}
\tilde{G}=
\begin{pmatrix}
g_{0,0} & g_{1,0} & g_{2,0} & \cdots & g_{m-1,0} & g_{0,1} & g_{1,1} & g_{2,1} & \cdots & g_{m-1,n-1}\\
0 & g_{1,0} & g_{2,0} & \cdots & g_{m-1,0} & 0 & g_{1,1} & g_{2,1} & \cdots & g_{m-1,n-1}\\
0 & 0 & g_{2,0} & \cdots & g_{m-1,0} & 0 & 0 & g_{2,1} & \cdots  & g_{m-1,n-1}\\
\vdots & \vdots & \vdots & \ddots & \vdots & \vdots & \vdots & \vdots & \ddots &  \vdots\\
0 & 0 & \cdots & 0 & g_{m-1,0} & 0 & 0 & \cdots & 0 & g_{m-1,n-1}\\
0 & 0 & 0 & ...& 0 & g_{0,1} & g_{1,1} & g_{2,1} & ... & g_{m-1,n-1}\\
\vdots & \vdots & \vdots & \ddots & \vdots & \vdots & \vdots & \vdots & \ddots &  \vdots\\
0 & 0 & 0 & ...& 0 & 0 & 0 & 0 & ... & g_{m-1,n-1}\\
\end{pmatrix}.
\end{equation*}
Denote $D=diag(g(0,0), g(1,0), \dots, g(m-1,n-1))$ and $A$ the upper triangular matrix of ones with dimensions $m\times m$, i.e.
\begin{equation*}
A= 
\begin{pmatrix}
1&1&1&\cdots&1\\
0&1&1&\cdots&1\\
0&0&1&\cdots&1\\
\vdots&\vdots&\vdots&\ddots&\vdots\\
0&0&0&\cdots&1 
\end{pmatrix}_{m\times m}.
\end{equation*}
Then it holds
\begin{equation*}
\tilde{G}= 
\begin{pmatrix}
A&A&A&\cdots&A\\
0&A&A&\cdots&A\\
0&0&A&\cdots&A\\
\vdots&\vdots&\vdots&\ddots&\vdots\\
0&0&0&\cdots&A 
\end{pmatrix}
D.
\end{equation*}
Matrix $\tilde{G}$ is an upper triangular matrix with dimensions $nm\times nm$, in its each column there is only one particular $g(i,j)$ (several times). Its diagonal consists of elements $g(i,j)$, each of them occurring just once. 

The notation established above allows us to write
\begin{equation*}
v\left(\mathbb{P}^{(1)}\right)=\frac{\tilde{G}v\left(\mathbb{P}^*\right)}{\lVert\tilde{G}v\left(\mathbb{P}^*\right)\rVert_1},   
\end{equation*}
and the $k$-th descendant can be expressed in its vector form as
\begin{equation*}
v\left(\mathbb{P}^{(k)}\right)=\frac{\tilde{G}v\left(\mathbb{P}^{(k-1)}\right)}{\lVert\tilde{G}v\left(\mathbb{P}^{(k-1)}\right)\rVert_1}=\frac{\tilde{G}^kv\left(\mathbb{P}^{*}\right)}{\lVert\tilde{G}^kv\left(\mathbb{P}^{*}\right)\rVert_1}.     
\end{equation*}

If the assumptions under which the power method converges are satisfied (i.e. a finite support of the parent distribution, a unique dominant eigenvalue of a diagonalizable matrix $\tilde{G}$, a suitable starting vector $\mathbb{P}^*$), the sequence
\begin{equation*}
\frac{\tilde{G}v\left(\mathbb{P}^*\right)}{\lVert\tilde{G}v\left(\mathbb{P}^*\right)\rVert_2}, \frac{\tilde{G}v\left(\mathbb{P}^{(1)}\right)}{\lVert\tilde{G}v\left(\mathbb{P}^{(1)}\right)\rVert_2}, ..., \frac{\tilde{G}v\left(\mathbb{P}^{(k)}\right)}{\lVert\tilde{G}v\left(\mathbb{P}^{(k)}\right)\rVert_2}, ...
\end{equation*}
converges to the unit dominant eigenvector of matrix $\tilde{G}$.
The eigenvalues of the upper triangular matrix are its diagonal elements, so the dominant eigenvalue is unique if and only if the greatest absolute value of $g(i,j)$ is unique.
We will obtain the limit distribution by multiplying the dominant unit vector by a normalization constant, i.e. the limit distribution will be 
\begin{equation*}
v(\mathbb{P}^{(\infty)})=\lim_{k\rightarrow\infty}\frac{v(\mathbb{P}^{(k)})}{\lVert v(\mathbb{P}^{(k)})\rVert_1}=\lim_{k\rightarrow\infty}\frac{\tilde{G}^kv(\mathbb{P}^{*})}{\lVert (\tilde{G}^kv(\mathbb{P}^{*}))\rVert_1}.
\end{equation*}

%

\vskip 3mm

\noindent 4.  EXAMPLES

\vskip 3mm
\noindent 4.1.  LIMIT DISTRIBUTION
  
 Let $N_1, N_2, N_3 \in \mathbb{N}$, $k\in \{1, 2, ..., N_3\}$ and $N=N_1+N_2+N_3$. Vector ${X\choose Y}$ has a bivariate inverse hypergeometric distribution (see Johnson et al. (1997)) if 
\begin{equation*}\label{BHD}
P(X=x,Y=y)=\frac{N_3-k+1}{N-(x+y+k-1)}\frac{\displaystyle{N_1\choose x}{N_2\choose y}{N_3\choose k-1}}{\displaystyle{N\choose {x+y+k-1}}},
\end{equation*}
for $x=0,1, 2, \dots, N_1$, $y=0,1, 2, \dots, N_2$.
 
 We will consider such a function $g(i,j)$ which leaves the bivariate inverse hypergeometric distribution unchanged. We choose parameter values $N_1=N_2=2$, $N_3=5$, $k=2$, i.e.

\begin{equation*}
\mathbb{P}=
\begin{pmatrix}
\frac{5}{18} &\frac{10}{63} & \frac{5}{126}\\
\frac{10}{63} &\frac{10}{63} &\frac{4}{63}\\
\frac{5}{126} &\frac{4}{63} &\frac{5}{126}
\end{pmatrix}.
\end{equation*}


The corresponding matrix $\tilde{G}$ is

\begin{equation*}
\tilde{G}=
\begin{pmatrix}
\frac{3}{7} & \frac{3}{20}& -\frac{3}{5}& \frac{3}{20} & \frac{9}{20} & \frac{3}{8}&-\frac{3}{5} & \frac{3}{8} & 1	\\	
0 & \frac{3}{20}& -\frac{3}{5}& 0 & \frac{9}{20} & \frac{3}{8}&0 & \frac{3}{8} & 1	\\
0 & 0& -\frac{3}{5}& 0 & 0 & \frac{3}{8}&0 & 0 & 1	\\
0 & 0& 0& \frac{3}{20} & \frac{9}{20} & \frac{3}{8}&-\frac{3}{5} & \frac{3}{8} & 1	\\
0 & 0& 0& 0 & \frac{9}{20} & \frac{3}{8}&0 & \frac{3}{8} & 1	\\
0 & 0& 0& 0 & 0 & \frac{3}{8}&0 & 0 & 1	\\
0 & 0& 0& 0 & 0 & 0&-\frac{3}{5} & \frac{3}{8} & 1	\\
0 & 0& 0& 0 & 0 & 0&0 & \frac{3}{8} & 1	\\
0 & 0& 0& 0 & 0 & 0&0 & 0 & 1	\\
\end{pmatrix}.
\end{equation*}

The conditions that matrix $\tilde{G}$ must be diagonalizable and it must have the unique dominant eigenvalue are satisfied in this case.
If we start from any suitable probability vector (it can not be ortoghonal to the space of the dominant eigenvalue), the iterated partial summations will converge to the unit eigenvector corresponding to the dominant eigenvalue of the matrix $\tilde{G}$. 	
In this case is the dominant eigenvalue $1$, so a multiple of its corresponding eigenvector will be the limit distribution
\begin{equation*}
\mathbb{P}^{(\infty)}=
\begin{pmatrix}
\frac{5}{18} &\frac{10}{63} & \frac{5}{126}\\
\frac{10}{63} &\frac{10}{63} &\frac{4}{63}\\
\frac{5}{126} &\frac{4}{63} &\frac{5}{126}
\end{pmatrix}.
\end{equation*}
We remind that (almost, with the exception of vectors orthogonal to the space of the dominant eigenvalue) regardless of the parent distribution, the limit distribution $\mathbb{P}^{(\infty)}$ is the bivariate hypergeometric distribution with the parameters $N_1=N_2=2$, $N_3=5$, $k=2$ which in our example determines the function $g(i,j)$.

\vskip 3mm

\noindent 4.2. OSCILLATION

There are also sequences of descendant distributions which do not converge. Let the parent be the bivariate hypergeometric distribution (see Johnson et al. (1997)) with the parameters $N_1=N_2=1$, $N_3=2$, $n=1$, i.e. 
\begin{equation*}
\mathbb{P}^* =
\begin{pmatrix}
\frac{1}{2} & \frac{1}{4} \\ \frac{1}{4} & 0 
\end{pmatrix}
\end{equation*}
and let the matrix $\tilde{G}$ be
\begin{equation*}
\tilde{G} = 
\begin{pmatrix}
-1 & 1 & 1 & 0 \\
0 & 1 & 0 & 0 \\
0 & 0 & 1 & 0 \\
0 & 0 & 0 & 0 \\
\end{pmatrix}.
\end{equation*}
After the first partial summation we obtain
\begin{equation*}
\mathbb{P}^{(1)}=
\begin{pmatrix}
0 & \frac{1}{2} \\ \frac{1}{2} & 0 
\end{pmatrix}, 
\end{equation*}
and after the second summation
\begin{equation*}
\mathbb{P}^{(2)}=
\begin{pmatrix}
\frac{1}{2} & \frac{1}{4} \\ \frac{1}{4} & 0 
\end{pmatrix},
\end{equation*}
i.e. the distribution identical to the parent $\mathbb{P}^*$.
\vskip 3mm

\noindent BIBLIOGRAPHY
\vskip 3mm

\noindent Golub H. G. and Van Loan Ch. F. (1996). {\it Matrix Computations.} Baltimore, London: The Johns Hopkins University Press.

\vskip 3mm

\noindent Johnson, N.L., Kemp, A.W. and Kotz, S. (2005). \it{Univariate Discrete Distributions.} \rm Hoboken (NJ): Wiley.

\vskip 3mm

\noindent Johnson, N.L., Kotz, S. and Balakrishnan, N. (1997). {\it{D}iscrete {M}ultivariate {D}istributions}. Hoboken (NJ): Wiley.

\vskip 3mm

\noindent Ko\v{s}\v{c}ov\'{a} M., Harman R. and Ma\v{c}utek J. (2018). Iterated partial summations applied to finite-support discrete distributions.

\noindent
http://www.iam.fmph.uniba.sk/ospm/Harman/KoscovaHarmanMacutek2018preprint.pdf

\noindent
(accessed on 30-Dec-2018)

\vskip 3mm

\noindent Kotz, S. and Johnson, N. L. (1991). A note on renewal (partial sums) distributions for discrete variables. {\it Statistics \& Probability Letters,} {\bf 12}, 229--231.\rm

\vskip 3mm

\noindent  Ma\v{c}utek J. (2006). A limit property of the geometric distribution.
{\it Theory of Probability and its Applications},  {\bf  50(2)}, 316--319.\rm

\vskip 3mm

\noindent  Ma\v{c}utek, J. (2003). On two types of partial summations. {\it Tatra Mountains	Mathematical Publications,} {\bf 26}, 403--410.\rm
 
\vskip 3mm

\noindent Wimmer, G. and Altmann, G. (2000). On the generalization of the STER distribution applied to generalized hypergeometric parents. {\it Acta Universitatis Palackianae Olomucensis Facultas Rerum Naturalium Mathematica,} {\bf 39(1)}, 215--247.\rm

\vskip 3mm

\noindent Wimmer, G. and Ma\v{c}utek, J. (2012). New integrated view at partial-sums distributions. {\it Tatra Mountains	Mathematical Publications,} {\bf 51}, 183--190.\rm


\noindent  
\vskip 1in

\end{document}